\documentclass[12pt,oneside]{amsart}

\usepackage{setspace}

\usepackage{amssymb}   

\usepackage{amsmath}

\usepackage{mathrsfs}

\usepackage{stmaryrd}

\usepackage{amsthm}
\usepackage{newlfont}
\usepackage{amscd}
\usepackage{mathtools}
\usepackage{bm}
\usepackage{tikz-cd}
\usepackage{graphicx}

\usepackage{hyperref}

\usepackage{enumitem}

\usepackage[all]{xy}

\usepackage{textcomp}

\usepackage{amsbsy}

\newtheorem{thm}{Theorem}[section]

\newtheorem{thm-defn}[thm]{Theorem/Definition}
\newtheorem{lem}[thm]{Lemma}
\newtheorem{prop}[thm]{Proposition}

\theoremstyle{definition}
\newtheorem{defn}[thm]{Definition}

\newtheorem{ques}[thm]{Question}

\theoremstyle{remark}
\newtheorem{rem}[thm]{Remark}

\usepackage{relsize}
\usepackage[bbgreekl]{mathbbol}
\usepackage{amsfonts}

\DeclareSymbolFontAlphabet{\mathbb}{AMSb}
\DeclareSymbolFontAlphabet{\mathbbl}{bbold}
\newcommand{\Prism}{{\mathlarger{\mathbbl{\Delta}}}}

\numberwithin{equation}{section}

\begin{document}

\pagenumbering{arabic}

\title{A note on purity of crystalline local systems}
\author{Yong Suk Moon}
\address{Beijing Institute of Mathematical Sciences and Applications, Beijing 101408, China}
\email{ysmoon@bimsa.cn}

\begin{abstract}
In this short note, we prove a purity result for crystalline local systems on a smooth $p$-adic affine formal scheme. Our method is based on the prismatic description of crystalline local systems \cite{du-liu-moon-shimizu-completed-prismatic-F-crystal-loc-system} (cf. \cite{guo-reinecke-prism-F-crys}).    
\end{abstract}

\maketitle

\section{Introduction} \label{sec:intro}

Let $K$ be a complete discrete valued field of mixed characteristic $(0, p)$ with the ring of integers $\mathcal{O}_K$ and perfect residue field $k$. Denote $K_0 = W(k)[p^{-1}]$ and $G_K = \mathrm{Gal}(\overline{K}/K)$ where $\overline{K}$ is an algebraic closure of $K$. To any finite dimensional continuous $\mathbf{Q}_p$-representation $V$ of $G_K$, Fontaine attached a $K_0$-vector space $D_{\mathrm{cris}}(V)$. More precisely, Fontaine introduced the crystalline period ring $\mathbf{B}_{\mathrm{cris}}$ equipped with a natural Frobenius and $G_K$-action, and considered $D_{\mathrm{cris}}(V) \coloneqq (V\otimes_{\mathbf{Q}_p} \mathbf{B}_{\mathrm{cris}})^{G_K}$ (\cite{fontaine-barsotti-tate}, \cite{fontaine-p-adic-periods}). We have $\dim_{K_0} D_{\mathrm{cris}}(V) \leq \dim_{\mathbf{Q}_p} V$ in general, and $V$ is called \textit{crystalline} if the equality holds. The underlying motivation of this notion comes from its close connection to good reduction. If $X$ is a proper smooth scheme over $K$ with good reduction, i.e. if there is a proper smooth scheme $\mathcal{X} / \mathcal{O}_K$ such that $\mathcal{X}\times_{\mathcal{O}_K} K = X$, then the $G_K$-representation $H^i_{\text{\'et}}(X_{\overline{K}}, \mathbf{Q}_p)$ is crystalline and $D_{\mathrm{cris}}(H^i_{\text{\'et}}(X_{\overline{K}}, \mathbf{Q}_p)) \cong H^i_{\mathrm{cris}}(\mathcal{X}_k / W(k))[p^{-1}]$. 

When $X$ is an abelian variety over $K$, the converse also holds: $X$ has good reduction if and only if $H^1_{\text{\'et}}(X_{\overline{K}}, \mathbf{Q}_p)$ is crystalline (\cite{coleman-iovita-frobenius-monodromy-abelian-variety}, \cite{mokrane-hyodo-kato-cohomology}). However, for abelian schemes over regular bases of mixed characteristic, an interesting discrepancy occurs in terms of purity. Motivated by Grothendieck's work on Nagata--Zariski purity (\cite{grothendieck-sga2}), we can ask the following question.

\begin{ques}
Let $R = \mathcal{O}_K[\![t_1, \ldots, t_d]\!]$ (with $d \geq 1$), and let $\mathfrak{m} \subset R$ be the maximal ideal. Given any abelian scheme over $\mathrm{Spec} R ~\backslash ~\{\mathfrak{m}\}$, does it extend uniquely to an abelian scheme over $\mathrm{Spec} R$? 	
\end{ques}

\noindent The answer to this question is positive if the ramification index $e = [K : K_0] \leq p-1$, but is negative if $e \geq p$ (\cite{vasiu-zink-purity}). We remark that when $e \leq p-1$, the purity result implies the uniqueness of integral canonical models of Shimura varieties (\cite[Cor.~30]{vasiu-zink-purity}).

On the other hand, for \textit{arbitrary} ramification index $e$, one can still study an analogous question on purity for geometric families of crystalline representations as follows. Faltings introduced the notion of crystalline $p$-adic \'etale local systems on the generic fiber of a smooth proper scheme over $\mathcal{O}_K$ (\cite{faltings}). Furthermore, for certain affine schemes which include the cases we study in this paper, Brinon studied the foundation for crystalline local systems \textit{\`a la} Fontaine by generalizing the construction of the crystalline period ring $\mathbf{B}_{\mathrm{cris}}$ (\cite{brinon-crys-rep-imperfect-residue}, \cite{brinon-relative}). We consider the small affine case, i.e. when $R$ is the $p$-adic completion of an \'etale algebra over $\mathcal{O}_K[T_1^{\pm 1}, \ldots, T_d^{\pm 1}]$ such that $\mathrm{Spec}(R/\pi R)$ is connected. Denote $\mathcal{G}_R \coloneqq \pi_1^{\text{\'et}}(\mathrm{Spec} R[p^{-1}])$. 

Let $\pi \in \mathcal{O}_K$ be a uniformizer, and let $\mathcal{O}_L$ be the $p$--adic completion of $R_{(\pi)}$. Denote by $\mathcal{G}_{O_L}$ the absolute Galois group of $L = \mathcal{O}_L[p^{-1}]$. Choose a geometric point of $\mathrm{Spec} (L)$, which gives a geometric point of $\mathrm{Spec}(R[p^{-1}])$ via the map $\mathrm{Spec} (L) \rightarrow \mathrm{Spec}(R[p^{-1}])$. By the change of paths for \'etale fundamental groups, we then have a continuous map of Galois groups $\mathcal{G}_{\mathcal{O}_L} \rightarrow \mathcal{G}_R$. For a finite dimensional continuous $\mathbf{Q}_p$-representation $V$ of $\mathcal{G}_R$, we refer the reader to \cite[\S 8.2]{brinon-relative} for the definition of $V$ being Hodge--Tate, de Rham, or crystalline. In this paper, we prove the the following purity statement.

\begin{thm} \label{thm:main-intro}
Let $R$ be the $p$-adic completion of an \'etale algebra over $\mathcal{O}_K[T_1^{\pm 1}, \ldots, T_d^{\pm 1}]$, and let $V$ be a finite dimensional continuous $\mathbf{Q}_p$-representation of $\mathcal{G}_R$. Then $V$ is crystalline if and only if $V|_{\mathcal{G}_{\mathcal{O}_L}}$ is crystalline.	
\end{thm}

\noindent Note that the ``only if'' part of the above theorem follows directly from the definition.

\begin{rem}
In \cite[Thm.~5.4.8]{tsuji-cryst-shvs}, Tsuji has already proved that if $V$ is de Rham and $V|_{\mathcal{G}_{\mathcal{O}_L}}$ is crystalline, then $V$ is crystalline. Furthermore, the purity of de Rham representations is expected to hold; it is expected that a similar argument as in the proof of \cite[Thm.~1.5 (ii)]{liu-zhu-rigidity} would imply that $V$ is de Rham if and only if $V|_{\mathcal{G}_{\mathcal{O}_L}}$ is de Rham. Another possible approach is based on the method in \cite{tsuji-Hodge-Tate-purity}. By \cite[Thm.~9.1]{tsuji-Hodge-Tate-purity}, $V$ is Hodge--Tate if and only if $V|_{\mathcal{G}_{\mathcal{O}_L}}$ is Hodge--Tate. It is expected that a similar argument as in \textit{loc. cit.} can be used via the results in \cite{andreatta-Brinon} to show $V$ is de Rham if and only if $V|_{\mathcal{G}_{\mathcal{O}_L}}$ is de Rham. Since any crystalline representation is de Rham, the expected purity of de Rham representations combined with the result of Tsuji \cite[Thm.~5.4.8]{tsuji-cryst-shvs} would imply Theorem~\ref{thm:main-intro}.
\end{rem}
  
Our method in this paper is completely different from the ones in the above remark. We employ the prismatic description of crystalline local systems given in \cite{du-liu-moon-shimizu-completed-prismatic-F-crystal-loc-system} (cf. \cite{guo-reinecke-prism-F-crys}).\\

\noindent \textbf{Notation.} Fix a prime $p$. Let $k$ be a perfect field of characteristic $p$, and let $K$ be a finite totally ramified extension of $K_0 \coloneqq W(k)[p^{-1}]$ with ring of integers $\mathcal{O}_K$. Fix a uniformizer $\pi \in \mathcal{O}_K$, and let $E = E(u) \in W(k)[u]$ be the monic minimal polynomial of $\pi$. 

For a ring $A$ and a finitely generated ideal $J \subset A$, the $J$-\textit{adic completion} of an $A$-module means the classical completion. Similarly, being $J$-\textit{adically complete} or $J$-\textit{complete} is in the classical sense. For a $\mathbf{Z}$-module $M$, we denote its $p$-adic completion by $M^{\wedge}_p$. 

A $p$-\textit{adically completed \'etale map} from a $p$-adically complete ring $B$ refers to the $p$-adic completion of an \'etale map from $B$. Write $\mathcal{O}_K\langle T_1^{\pm 1}, \ldots, T_d^{\pm 1}\rangle$ for the $p$-adic completion of the Laurent polynomial ring $\mathcal{O}_K[T_1^{\pm 1}, \ldots, T_d^{\pm 1}]$, and similarly for $W(k)\langle T_1^{\pm 1}, \ldots, T_d^{\pm 1}\rangle$.

For an element $a$ of a $\mathbf{Q}$-algebra $A$ and $n \geq 0$, write $\gamma_n(a)$ for the element $\frac{a^n}{n!} \in A$.

\section*{Acknowledgements}

I would like to thank Heng Du, Tong Liu, and Koji Shimizu for helpful discussions. I also thank the anonymous referees for many valuable suggestions to improve the paper.

\section{Absolute prismatic site \& Kisin descent datum} \label{sec:Kisin-descent-datum}

We first recall some of the main results in \cite{du-liu-moon-shimizu-completed-prismatic-F-crystal-loc-system} on the equivalence between the category of crystalline local systems and the category of \textit{Kisin descent data}. As in \cite[Assumption~2.9]{du-liu-moon-shimizu-completed-prismatic-F-crystal-loc-system}, we consider the cases when the base ring is either small over $\mathcal{O}_K$ or a complete discrete valuation ring. 

Let $R$ be a $p$-adically completed \'etale algebra over $\mathcal{O}_K\langle T_1^{\pm 1}, \ldots, T_d^{\pm 1}\rangle$ for some $d \geq 0$ such that $\mathrm{Spec}(R/\pi R)$ is connected. There exists a subring $R_0 \subset R$ such that $R_0$ is $p$-adically completed \'etale over $W(k)\langle T_1^{\pm 1}, \ldots, T_d^{\pm 1}\rangle$ and $R_0\otimes_{W(k)} \mathcal{O}_K = R$ (see e.g. \cite[Lem.~2.9]{guo-reinecke-prism-F-crys}). Let $\varphi\colon R_0 \rightarrow R_0$ be the (unique) lift of Frobenius on $R_0/pR_0$ with $\varphi(T_i) = T_i^p$. Let $\mathcal{O}_{L_0}$ be the $p$-adic completion of $(R_0)_{(p)}$ equipped with the Frobenius induced from $\varphi$ on $R_0$. Note that $\mathcal{O}_{L_0}$ is a complete discrete valuation ring whose residue field has a finite $p$-basis given by $\{T_1, \ldots, T_d\}$. We have a natural injective map $R_0 \rightarrow \mathcal{O}_{L_0}$ compatible with $\varphi$. This extends $\mathcal{O}_K$-linearly to $R \rightarrow \mathcal{O}_L \coloneqq \mathcal{O}_{L_0}\otimes_{W(k)} \mathcal{O}_K$.

\medskip 

\noindent\textbf{Assumption.}
In the following, we will assume the base ring $S$ is either $R$ or $\mathcal{O}_L$. Denote $S_0 = R_0$ (resp. $S_0 = \mathcal{O}_{L_0}$) when $S = R$ (resp. $S = \mathcal{O}_L$).	

\begin{defn}[{\cite[Def.~3.2]{bhatt-scholze-prismaticcohom-v3}}]
A \textit{bounded prism} is a pair $(A, I)$ where $A$ is a $\delta$-ring (cf. \cite[Def.~2.1]{bhatt-scholze-prismaticcohom-v3}) and $I \subset A$ is an invertible ideal such that $p \in I+\varphi(I) A$, $A/I$ has bounded $p^{\infty}$-torsion, and $A$ is $(p, I)$-complete (see \cite[Lem.~3.7 (1)]{bhatt-scholze-prismaticcohom-v3}). Here, $\varphi\colon A \rightarrow A$ is given by $\varphi(x) = x^p+p\delta(x)$.	
\end{defn}

\begin{defn}[{\cite[Def.~2.3]{bhatt-scholze-prismaticFcrystal}}]
The \textit{absolute prismatic site} $S_{\Prism}$ of the $p$-adic formal scheme $\mathrm{Spf} S$ consists of the pairs $((A, I), \mathrm{Spf} A/I \rightarrow \mathrm{Spf} S)$ where $(A, I)$ is a bounded prism and $\mathrm{Spf} A/I \rightarrow \mathrm{Spf} S$ is a morphism of $p$-adic formal schemes. For simplicity, we often omit the structure map $\mathrm{Spf} A/I \rightarrow \mathrm{Spf} S$ and write $(A, I) \in S_{\Prism}$. The morphisms are the opposite of morphisms of bounded prisms compatible with the structure maps to $\mathrm{Spf} S$. We equip $S_{\Prism}$ with the topology given by $(p, I)$-completely faithfully flat morphisms of bounded prisms $(A, I) \rightarrow (B, J)$.   	
\end{defn}

An important object in $S_{\Prism}$ is the \textit{Breuil--Kisin prism} given as follows. Denote $\mathfrak{S}_S = S_0[\![u]\!]$ equipped with the Frobenius extending that on $S_0$ such that $\varphi(u) = u^p$. Then $(\mathfrak{S}_S, (E))$ is a bounded prism, and it is an object in $S_{\Prism}$ via $\mathfrak{S}_S/(E) \cong S$. 

The self-product of $(\mathfrak{S}_S, (E))$ in $S_{\Prism}$ exists and is computed in \cite[Ex.~3.4]{du-liu-moon-shimizu-completed-prismatic-F-crystal-loc-system}, which we briefly explain. Write $\mathfrak{S}_S\widehat{\otimes}_{\mathbf{Z}_p}\mathfrak{S}_S$ for the $p$-complete tensor product  equipped with the induced $\otimes$-product Frobenius, and consider $d\colon \mathfrak{S}_S\widehat{\otimes}_{\mathbf{Z}_p}\mathfrak{S}_S \rightarrow S$ given by the composite $\mathfrak{S}_S\widehat{\otimes}_{\mathbf{Z}_p}\mathfrak{S}_S \rightarrow \mathfrak{S}_S \rightarrow \mathfrak{S}_S/(E) \cong S$ where the first map is the multiplication. Let $J$ be the kernel of $d$, and consider
\[
\mathfrak{S}_S^{(1)} \coloneqq (\mathfrak{S}_S\widehat{\otimes}_{\mathbf{Z}_p} \mathfrak{S}_S)\biggl \{\frac{J}{E}\biggr\}_\delta^{\wedge}.
\]
Here the $\mathfrak{S}_S$-algebra structure of $\mathfrak{S}_S^{(1)}$ is given by $a \mapsto a\otimes 1$, and $\{\cdot \}^{\wedge}_{\delta}$ means adjoining elements in the category of derived $(p, E)$-complete simplicial $\delta$-$\mathfrak{S}_S$-algebras. Note that $E$ in $\{\frac{J}{E}\}^{\wedge}_{\delta}$ denotes $E\otimes 1$, but using $1\otimes E$ yields the same $\mathfrak{S}_S^{(1)}$. We have $(\mathfrak{S}_S^{(1)}, (E)) \in S_{\Prism}$, and it is the self-product of $(\mathfrak{S}_S, (E))$ in $S_{\Prism}$. Similarly, the self-triple-product $(\mathfrak{S}_S^{(2)}, (E))$ of the Breuil--Kisin prism exists in $S_{\Prism}$. Write $p_i\colon \mathfrak{S}_S \rightarrow \mathfrak{S}_S^{(1)}$ with $i = 1, 2$ and $q_i\colon \mathfrak{S}_S \rightarrow \mathfrak{S}_S^{(2)}$ with $i = 1, 2, 3$ for the projection maps. Note that by the rigidity of maps of prisms (\cite[Lem.~3.5]{bhatt-scholze-prismaticcohom-v3}), we have $(p_1(E)) = (E) = (p_2(E))$ as ideals of $\mathfrak{S}_S^{(1)}$, and similarly $(q_i(E)) = (E)$ as ideals of $\mathfrak{S}^{(2)}$ for each $i = 1, 2, 3$. 

The Breuil--Kisin prism covers the final object of $\mathrm{Shv}(S_{\Prism})$, and thus a crystal on $S_{\Prism}$ can be described by a $\mathfrak{S}_S$-module with a descent datum involving the self-product and self-triple-product. For example, \textit{completed prismatic} $F$-\textit{crystals} on $S_{\Prism}$ given in \cite[Def.~3.16]{du-liu-moon-shimizu-completed-prismatic-F-crystal-loc-system}) can be described in terms of \textit{Kisin descent data} defined below (\cite[Prop.~3.26]{du-liu-moon-shimizu-completed-prismatic-F-crystal-loc-system}).

\begin{defn}[cf. {\cite[Def.~3.14]{du-liu-moon-shimizu-completed-prismatic-F-crystal-loc-system}}] \hfill
\begin{itemize}
\item We say that a finite $\mathfrak{S}_S$-module $N$ is \textit{projective away from} $(p, E)$ if $N$ is $p$-torsion free, $N[p^{-1}]$ is projective over $\mathfrak{S}_S[p^{-1}]$, and $N[E^{-1}]^{\wedge}_p$ is projective over $\mathfrak{S}_S[E^{-1}]^{\wedge}_p$. 
\item We say a finite $\mathfrak{S}_S$-module $N$ is \textit{saturated} if $N$ is torsion free and $N = N[p^{-1}] \cap N[E^{-1}]$.
\item Let $N$ be a $\mathfrak{S}_S$-module equipped with a $\varphi$-semi-linear endomorphism $\varphi_N\colon N \rightarrow N$. We say $(N, \varphi_N)$ has \textit{finite} $E$-\textit{height} if $1\otimes\varphi_N \colon \mathfrak{S}_S\otimes_{\varphi, \mathfrak{S}_S}N \rightarrow N$ is injective and its cokernel is killed by a power of $E$. 	
\end{itemize}	
\end{defn}

\begin{rem}
When $S = \mathcal{O}_L$, any finite $\mathfrak{S}_{\mathcal{O}_L}$-module which is projective away from $(p, E)$ and saturated is free over $\mathfrak{S}_{\mathcal{O}_L}$, since $\mathfrak{S}_{\mathcal{O}_L}$ is a regular local ring of dimension $2$ (cf. \cite[Rem.~3.18]{du-liu-moon-shimizu-completed-prismatic-F-crystal-loc-system}).
\end{rem}

\begin{defn}[{\cite[Def.~3.25]{du-liu-moon-shimizu-completed-prismatic-F-crystal-loc-system}}] Let $\mathrm{DD}_{\mathfrak{S}_S}$ denote the category consisting of triples $(\mathfrak{M}, \varphi_{\mathfrak{M}}, f)$ (called \textit{Kisin descent datum}) where
\begin{itemize}
\item $\mathfrak{M}$ is a finite $\mathfrak{S}_S$-module that is projective away from $(p, E)$ and saturated;
\item $\varphi_{\mathfrak{M}}\colon \mathfrak{M} \rightarrow \mathfrak{M}$ is a $\varphi$-semi-linear endomorphism such that $(\mathfrak{M}, \varphi_{\mathfrak{M}})$ has finite $E$-height;
\item $f\colon \mathfrak{S}_S^{(1)}\otimes_{p_1, \mathfrak{S}_S} \mathfrak{M} \stackrel{\cong}{\rightarrow} \mathfrak{S}_S^{(1)}\otimes_{p_2, \mathfrak{S}_S} \mathfrak{M}$ is an isomorphism of $\mathfrak{S}_S^{(1)}$-modules compatible with Frobenii and satisfies the cocycle condition over $\mathfrak{S}_S^{(2)}$ (i.e. if $p_{12}, p_{23}, p_{13}\colon \mathfrak{S}_S^{(1)} \rightarrow \mathfrak{S}_S^{(2)}$ denote the projections, then $p_{23}^*f \circ p_{12}^*f = p_{13}^* f$). 	
\end{itemize}	
\end{defn}

The main input we will need is the following theorem proved in \cite{du-liu-moon-shimizu-completed-prismatic-F-crystal-loc-system}. Recall that if a finite dimensional continuous $\mathbf{Q}_p$-representation $V$ of $\mathcal{G}_S\coloneqq \pi_1^{\text{\'et}}(\mathrm{Spec} S[p^{-1}])$ is crystalline, then it is de Rham, and we can attach a $S[p^{-1}]$-module $D_{\mathrm{dR}}(V)$ projective of rank equal to $\dim_{\mathbf{Q}_p} V$ (\cite[\S 8]{brinon-relative}). Then $D_{\mathrm{dR}}(V)$ is equipped with a decreasing exhaustive filtration by $S[p^{-1}]$-submodules $\mathrm{Fil}^i D_{\mathrm{dR}}(V)$, and the Hodge--Tate weights of $V$ are defined to be the integers $i$ such that $\mathrm{Fil}^i D_{\mathrm{dR}}(V) \neq \mathrm{Fil}^{i+1} D_{\mathrm{dR}}(V)$.   

\begin{thm}[{\cite[Prop.~3.26, Thm.~3.29]{du-liu-moon-shimizu-completed-prismatic-F-crystal-loc-system}}] \label{thm:equiv-crys-loc-syst-kisin-descent-data}
	The category $\mathrm{DD}_{\mathfrak{S}_S}$ is naturally equivalent to the category of $\mathbf{Z}_p$-lattices of crystalline representations of $\mathcal{G}_S$ with non-negative Hodge--Tate weights.
\end{thm}

\section{Purity of crystalline local systems} \label{sec:purity}

We prove Theorem \ref{thm:main-intro} in this section. Let $\mathcal{G}_{\mathcal{O}_L} \rightarrow \mathcal{G}_R$ be a map of Galois groups as in \S\ref{sec:intro}. Let $V$ be a finite dimensional continuous $\mathbf{Q}_p$--representation of $\mathcal{G}_R$ such that $V|_{\mathcal{G}_{\mathcal{O}_L}}$ is crystalline. By applying a suitable power of Tate twist (i.e. twist by a power of the $p$-adic cyclotomic character), we may assume that the Hodge--Tate weights of $V|_{\mathcal{G}_{\mathcal{O}_L}}$ are non-negative. Let $T \subset V$ be a $\mathcal{G}_R$-stable $\mathbf{Z}_p$-lattice. By \cite[Cor.~3.8]{bhatt-scholze-prismaticFcrystal}, we can naturally associate to $T$ an \'etale $\varphi$--module $\mathcal{M}$ which is finite projective over $\mathfrak{S}_R[E^{-1}]^{\wedge}_p$ together with a $\mathfrak{S}_R^{(1)}[E^{-1}]^{\wedge}_p$-linear isomorphism
\[
f_{\text{\'et}}\colon \mathfrak{S}_R^{(1)}[E^{-1}]^{\wedge}_p\otimes_{p_1, \mathfrak{S}_R[E^{-1}]^{\wedge}_p} \mathcal{M} \stackrel{\cong}{\rightarrow} \mathfrak{S}_R^{(1)}[E^{-1}]^{\wedge}_p\otimes_{p_2, \mathfrak{S}_R[E^{-1}]^{\wedge}_p} \mathcal{M},
\] 
which is compatible with $\varphi$ and satisfies the cocycle condition over $\mathfrak{S}_R^{(2)}[E^{-1}]^{\wedge}_p$. Here, we associate $\mathcal{M}$ and $T$ contravariantly following the convention in \cite{du-liu-moon-shimizu-completed-prismatic-F-crystal-loc-system} (see \cite[\S 3.4]{du-liu-moon-shimizu-completed-prismatic-F-crystal-loc-system}). 

Note that the map $R_0 \rightarrow \mathcal{O}_{L_0}$ extends to $\mathfrak{S}_R \rightarrow \mathfrak{S}_{\mathcal{O}_L}$ by $u \mapsto u$, which is compatible with Frobenius. Since $V|_{\mathcal{G}_{\mathcal{O}_L}}$ is crystalline with non-negative Hodge--Tate weights, by Theorem~\ref{thm:equiv-crys-loc-syst-kisin-descent-data} and \cite[Prop.~3.27]{du-liu-moon-shimizu-completed-prismatic-F-crystal-loc-system}, there exists a Kisin descent datum $(\mathfrak{M}_L, \varphi_{\mathfrak{M}_L}, f_L) \in \mathrm{DD}_{\mathfrak{S}_{\mathcal{O}_L}}$ over $\mathfrak{S}_{\mathcal{O}_L}$ such that we have a $\varphi$-compatible isomorphism $h\colon \mathfrak{M}_L\otimes_{\mathfrak{S}_{\mathcal{O}_L}} \mathfrak{S}_{\mathcal{O}_L}[E^{-1}]^{\wedge}_p \cong \mathcal{M}\otimes_{\mathfrak{S}_R[E^{-1}]^{\wedge}_p} \mathfrak{S}_{\mathcal{O}_L}[E^{-1}]^{\wedge}_p$ and the base change of the isomorphism
\[
f_L\colon \mathfrak{S}_{\mathcal{O}_L}^{(1)}\otimes_{p_1, \mathfrak{S}_{\mathcal{O}_L}} \mathfrak{M}_L \stackrel{\cong}{\rightarrow} \mathfrak{S}_{\mathcal{O}_L}^{(1)}\otimes_{p_2, \mathfrak{S}_{\mathcal{O}_L}} \mathfrak{M}_L
\]
to $\mathfrak{S}_{\mathcal{O}_L}^{(1)}[E^{-1}]^{\wedge}_p$ agrees with the base change of $f_{\text{\'et}}$ to $\mathfrak{S}_{\mathcal{O}_L}^{(1)}[E^{-1}]^{\wedge}_p$ (with respect to $h$). Let $\mathcal{M}_L = \mathfrak{M}_L\otimes_{\mathfrak{S}_{\mathcal{O}_L}} \mathfrak{S}_{\mathcal{O}_L}[E^{-1}]^{\wedge}_p$, and regard $\mathcal{M}$ as a $\mathfrak{S}_R[E^{-1}]^{\wedge}_p$-submodule of $\mathcal{M}_L$ via the isomorphism $h$.

Consider the $\mathfrak{S}$-module
\[
\mathfrak{M} \coloneqq \mathfrak{M}_L \cap \mathcal{M} \subset \mathcal{M}_L
\]  
equipped with the induced Frobenius $\varphi_{\mathfrak{M}}$. Note that $\mathfrak{M}$ is torsion free. By \cite[Prop.~4.20, 4.21]{du-liu-moon-shimizu-completed-prismatic-F-crystal-loc-system}, $\mathfrak{M}$ is finite over $\mathfrak{S}_R$, saturated, and has finite $E$-height with respect to $\varphi_{\mathfrak{M}}$. Furthermore, by \cite[Prop.~4.13, 4.26]{du-liu-moon-shimizu-completed-prismatic-F-crystal-loc-system}, $\mathfrak{M}$ is projective away from $(p, E)$ and we have natural $\varphi$-equivariant isomorphisms
\[
\mathfrak{M}\otimes_{\mathfrak{S}_R} \mathfrak{S}_{\mathcal{O}_L} \cong \mathfrak{M}_L ~\text{ and } ~\mathfrak{M}\otimes_{\mathfrak{S}_R} \mathfrak{S}_R[E^{-1}]^{\wedge}_p \cong \mathcal{M}.
\]
We claim that $f_{\text{\'et}}$ and $f_L$ induce an isomorphism
\[
f\colon \mathfrak{S}_R^{(1)}\otimes_{p_1, \mathfrak{S}_R} \mathfrak{M} \stackrel{\cong}{\rightarrow} \mathfrak{S}_R^{(1)}\otimes_{p_2, \mathfrak{S}_R} \mathfrak{M}
\]
so that $f$ is compatible with $f_{\text{\'et}}$ and $f_L$. For this, we need some preliminary facts.

\begin{lem} \label{lem:frakS(2)-to-frakSL(2)-injective}
The natural map 
\[
\mathfrak{S}_R^{(1)}/(p, E) \rightarrow \mathfrak{S}_{\mathcal{O}_L}^{(1)}/(p, E)
\]
is injective.	
\end{lem}

\begin{proof}
By \cite[Prop.~2.2.8 (2), 4.1.3]{du-liu-prismaticphiGhatmodule}, we have
\[
\mathfrak{S}_R^{(1)}/(E) \cong R[\gamma_i(z_j), i \geq 0, j = 0, \ldots, d]^{\wedge}_p
\]
and
\[
\mathfrak{S}_{\mathcal{O}_L}^{(1)}/(E) \cong \mathcal{O}_L[\gamma_i(z_j), i \geq 0, j = 0, \ldots, d]^{\wedge}_p
\]
where $z_0, \ldots, z_d$ can be considered as variables. Thus,
\[
\mathfrak{S}_R^{(1)}/(p, E) \cong R[\gamma_i(z_j), i \geq 0, j = 0, \ldots, d]/(p)
\]
and similarly for $\mathfrak{S}_{\mathcal{O}_L}^{(1)}/(p, E)$. Since $R/(p) \rightarrow \mathcal{O}_L/(p)$ is injective, the map 
\[
\mathfrak{S}_R^{(1)}/(p, E) \rightarrow \mathfrak{S}_{\mathcal{O}_L}^{(1)}/(p, E)
\]
is injective.
\end{proof}

\begin{lem}[{\cite[Cor.~3.6]{du-liu-moon-shimizu-completed-prismatic-F-crystal-loc-system}}] \label{lem:frakS(2)-basic-properties}
Let $S = R$ or $S = \mathcal{O}_L$. Then $\mathfrak{S}_S^{(1)}$ is $p$-torsion free and $E$-torsion free. Furthermore, 
\[
\mathfrak{S}_S^{(1)} = \mathfrak{S}_S^{(1)}[p^{-1}] \cap \mathfrak{S}_S^{(1)}[E^{-1}],
\]
and $\mathfrak{S}_S^{(1)}[E^{-1}]$ is $p$-adically separated. 
\end{lem}

\begin{lem} \label{lem:p-torsion-free}
Let $S = R$ or $S = \mathcal{O}_L$. Then $\mathfrak{S}_{S}^{(1)}[E^{-1}]^{\wedge}_p$ is $p$-torsion free.	
\end{lem}

\begin{proof}
By \cite[Lem.~3.5]{du-liu-moon-shimizu-completed-prismatic-F-crystal-loc-system}, $p_1\colon \mathfrak{S}_S \rightarrow \mathfrak{S}_S^{(1)}$ is classically faithfully flat. So the induced map $\mathfrak{S}_S[E^{-1}]^{\wedge}_p \rightarrow \mathfrak{S}_S^{(1)}[E^{-1}]^{\wedge}_p$	 is classically faithfully flat by \cite[Tag~0912]{stacks-project}. Since $\mathfrak{S}_S[E^{-1}]^{\wedge}_p$ is $p$-torsion free, the statement follows. 
\end{proof}

\begin{lem} \label{lem:E-complete}
Let $S = R$ or $S = \mathcal{O}_L$. Then $\mathfrak{S}_S^{(1)}/(p)$ is $E$-adically complete.	
\end{lem}

\begin{proof}
Note that $\mathfrak{S}_S^{(1)}$ is $(p, E)$-complete. Consider the exact sequence
\begin{equation} \label{eq:exact-seq}
0 \rightarrow p\mathfrak{S}_S^{(1)} \rightarrow \mathfrak{S}_S^{(1)} \rightarrow \mathfrak{S}_S^{(1)}/(p) \rightarrow 0.
\end{equation}
By Lemma~\ref{lem:frakS(2)-basic-properties}, $\mathfrak{S}_S^{(1)}$ is $p$-torsion free and $\mathfrak{S}_S^{(1)}/(p)$ is $E^n$-torsion free for each $n \geq 1$. In particular, the induced sequence 
\[
0 \rightarrow p\mathfrak{S}_S^{(1)}/E^n p\mathfrak{S}_S^{(1)} \rightarrow \mathfrak{S}_S^{(1)}/(E^n) \rightarrow \mathfrak{S}_S^{(1)}/(p, E^n) \rightarrow 0
\] 
is exact. By \cite[Tag~03CA]{stacks-project}, the sequence~(\ref{eq:exact-seq}) remains exact after $E$-completion. Since $p\mathfrak{S}_S^{(1)}$ is a free $\mathfrak{S}_S^{(1)}$-module of rank $1$, $p\mathfrak{S}_S^{(1)}$ is $E$-complete. Thus, $\mathfrak{S}_S^{(1)}/(p)$ is $E$-complete. 
\end{proof}

\begin{lem} \label{lem:injectivity-of-maps}
Let $S = R$ or $S = \mathcal{O}_L$. The natural maps
\[
\mathfrak{S}_S^{(1)}/(p) \rightarrow \mathfrak{S}_S^{(1)}[E^{-1}]/(p) ~~\text{and}~~\mathfrak{S}_S^{(1)} \rightarrow \mathfrak{S}_S^{(1)}[E^{-1}]^{\wedge}_p 
\]
are injective. Furthermore, the maps
\[
\mathfrak{S}_R^{(1)}/(p) \rightarrow \mathfrak{S}_{\mathcal{O}_L}^{(1)}/(p), ~~\mathfrak{S}_R^{(1)} \rightarrow \mathfrak{S}_{\mathcal{O}_L}^{(1)}, ~~\text{and}~~ \mathfrak{S}_R^{(1)}[E^{-1}]^{\wedge}_p \rightarrow \mathfrak{S}_{\mathcal{O}_L}^{(1)}[E^{-1}]^{\wedge}_p
\]
are injective.
\end{lem}

\begin{proof}
By Lemma~\ref{lem:frakS(2)-basic-properties}, $\{p, E\}$ form a regular sequence for $\mathfrak{S}_S^{(1)}$, and the maps $\mathfrak{S}_S^{(1)} \rightarrow \mathfrak{S}_S^{(1)}[E^{-1}]$ and $\mathfrak{S}_S^{(1)}[E^{-1}] \rightarrow \mathfrak{S}_S^{(1)}[E^{-1}]^{\wedge}_p$ are injective. Thus, the maps $\mathfrak{S}_S^{(1)}/(p) \rightarrow \mathfrak{S}_S^{(1)}[E^{-1}]/(p)$ and 	$\mathfrak{S}_S^{(1)} \rightarrow \mathfrak{S}_S^{(1)}[E^{-1}]^{\wedge}_p$ are injective.

Since $\mathfrak{S}_{\mathcal{O}_L}^{(1)}/(p)$ is $E$-torsion free, we deduce from Lemma~\ref{lem:frakS(2)-to-frakSL(2)-injective} inductively that the map $\mathfrak{S}_R^{(1)}/(p, E^n) \rightarrow \mathfrak{S}_{\mathcal{O}_L}^{(1)}/(p, E^n)$ is injective for each $n \geq 1$. By taking the inverse limit over $n$ giving the $E$-adic completions and using Lemma~\ref{lem:E-complete}, we have that the map $\mathfrak{S}_R^{(1)}/(p) \rightarrow \mathfrak{S}_{\mathcal{O}_L}^{(1)}/(p)$ is injective. Similarly, since $\mathfrak{S}_{\mathcal{O}_L}^{(1)}$ is $p$-torsion free and $\mathfrak{S}_R^{(1)}$ and $\mathfrak{S}_{\mathcal{O}_L}^{(1)}$ are $p$-complete, it follows that $\mathfrak{S}_R^{(1)} \rightarrow \mathfrak{S}_{\mathcal{O}_L}^{(1)}$ is injective. Furthermore, since $\mathfrak{S}_{\mathcal{O}_L}^{(1)}[E^{-1}]$ is $p$-torsion free and $\mathfrak{S}_R^{(1)}[E^{-1}]/(p) \rightarrow \mathfrak{S}_{\mathcal{O}_L}^{(1)}[E^{-1}]/(p)$ is injective, the map $\mathfrak{S}_R^{(1)}[E^{-1}]^{\wedge}_p \rightarrow \mathfrak{S}_{\mathcal{O}_L}^{(1)}[E^{-1}]^{\wedge}_p$ is injective.  
\end{proof}

\begin{prop} \label{prop:intersection-rings}
We have
\[
\mathfrak{S}_R^{(1)} = \mathfrak{S}_{\mathcal{O}_L}^{(1)} \cap \mathfrak{S}_R^{(1)}[E^{-1}]^{\wedge}_p
\]
as subrings of $\mathfrak{S}_{\mathcal{O}_L}^{(1)}[E^{-1}]^{\wedge}_p$.
\end{prop}

\begin{proof}
By Lemma~\ref{lem:injectivity-of-maps}, the map
\[
\mathfrak{S}_R^{(1)}/(p) \rightarrow (\mathfrak{S}_{\mathcal{O}_L}^{(1)} / (p)) \bigcap (\mathfrak{S}_R^{(1)}[E^{-1}] / (p)) 
\]
is injective, where the intersection is taken as subrings of $\mathfrak{S}_{\mathcal{O}_L}^{(1)}[E^{-1}]/(p)$. This map is also surjective by Lemma~\ref{lem:frakS(2)-to-frakSL(2)-injective}. Since $\mathfrak{S}_{\mathcal{O}_L}^{(1)}[E^{-1}]^{\wedge}_p$ is $p$-torsion free by Lemma~\ref{lem:p-torsion-free} and $\mathfrak{S}_R^{(1)}$ is $p$-complete, it follows that the map $\mathfrak{S}_R^{(1)} \rightarrow \mathfrak{S}_{\mathcal{O}_L}^{(1)} \cap \mathfrak{S}_R^{(1)}[E^{-1}]^{\wedge}_p$ is surjective. 
\end{proof}

Now, since $\mathfrak{M}[p^{-1}]$ is projective over $\mathfrak{S}_R[p^{-1}]$, we have
\[
(\mathfrak{S}_{\mathcal{O}_L}^{(1)}[p^{-1}]\otimes_{p_i, \mathfrak{S}_R[p^{-1}]} \mathfrak{M}[p^{-1}]) \bigcap (\mathfrak{S}_R^{(1)}[E^{-1}]^{\wedge}_p[p^{-1}]\otimes_{p_i, \mathfrak{S}_R[p^{-1}]} \mathfrak{M}[p^{-1}]) \cong (\mathfrak{S}_R^{(1)}[p^{-1}]\otimes_{p_i, \mathfrak{S}_R[p^{-1}]} \mathfrak{M}[p^{-1}])
\]
for $i = 1, 2$ by Proposition~\ref{prop:intersection-rings}. Thus, by \cite[Lem.~4.10]{du-liu-moon-shimizu-completed-prismatic-F-crystal-loc-system}, $f_{\text{\'et}}$ and $f_L$ induce a morphism
\[
f\colon \mathfrak{S}_R^{(1)}\otimes_{p_1, \mathfrak{S}_R} \mathfrak{M} \rightarrow \mathfrak{S}_R^{(1)}\otimes_{p_2, \mathfrak{S}_R} \mathfrak{M}. 
\]
Furthermore, since $f_{\text{\'et}}$ and $f_L$ are isomorphisms, it follows that $f$ obtained as their intersection is an isomorphism. Since $f_{\text{\'et}}$ is compatible with Frobenius, so is $f$. It remains to show that $f$ satisfies the cocycle condition over $\mathfrak{S}_R^{(2)}$.

\begin{lem} \label{lem:injectivity-for-frakS^(2)}
For each $i = 1, 2, 3$, the natural map
\[
\mathfrak{S}_R^{(2)}\otimes_{q_i, \mathfrak{S}_R} \mathfrak{M} \rightarrow \mathfrak{S}_R^{(2)}[E^{-1}]^{\wedge}_p\otimes_{q_i, \mathfrak{S}_R} \mathfrak{M}
\]	
is injective.
\end{lem}

\begin{proof}
First note that $q_i\colon \mathfrak{S}_R \rightarrow \mathfrak{S}_R^{(2)}$ is classically faithfully flat by \cite[Lem.~3.5]{du-liu-moon-shimizu-completed-prismatic-F-crystal-loc-system}. So by the same argument as in \cite[Cor.~3.6 Pf.]{du-liu-moon-shimizu-completed-prismatic-F-crystal-loc-system}, we deduce that $\mathfrak{S}_R^{(2)}$ is $p$-torsion free and $E$-torsion free, and $\mathfrak{S}_R^{(2)}[E^{-1}]$ is $p$-adically separated. In particular, the map $\mathfrak{S}_R^{(2)} \rightarrow \mathfrak{S}_R^{(2)}[E^{-1}]^{\wedge}_p$ is injective.

Furthermore, since $\mathfrak{M} \rightarrow \mathfrak{M}[p^{-1}]$ is injective, $\mathfrak{S}_R^{(2)}\otimes_{q_i, \mathfrak{S}_R} \mathfrak{M} \rightarrow \mathfrak{S}_R^{(2)}\otimes_{q_i, \mathfrak{S}_R} \mathfrak{M}[p^{-1}]$ is injective. The map $\mathfrak{S}_R^{(2)}\otimes_{q_i, \mathfrak{S}_R} \mathfrak{M}[p^{-1}] \rightarrow \mathfrak{S}_R^{(2)}[E^{-1}]^{\wedge}_p\otimes_{q_i, \mathfrak{S}_R} \mathfrak{M}[p^{-1}]$ is injective since $\mathfrak{M}[p^{-1}]$ is projective over $\mathfrak{S}_R[p^{-1}]$. Thus, the composite 
\[
\mathfrak{S}_R^{(2)}\otimes_{q_i, \mathfrak{S}_R} \mathfrak{M} \rightarrow \mathfrak{S}_R^{(2)}[E^{-1}]^{\wedge}_p\otimes_{q_i, \mathfrak{S}_R} \mathfrak{M}[p^{-1}]
\]  
is injective. Since this map factors through $\mathfrak{S}_R^{(2)}\otimes_{q_i, \mathfrak{S}_R} \mathfrak{M} \rightarrow \mathfrak{S}_R^{(2)}[E^{-1}]^{\wedge}_p\otimes_{q_i, \mathfrak{S}_R} \mathfrak{M}$, the statement follows.
\end{proof}

Since $f_{\text{\'et}}$ satisfies the cocycle condition over $\mathfrak{S}_R^{(2)}[E^{-1}]^{\wedge}_p$, we deduce from Lemma~\ref{lem:injectivity-for-frakS^(2)} that $f$ satisfies the cocycle condition over $\mathfrak{S}_R^{(2)}$. By Theorem~\ref{thm:equiv-crys-loc-syst-kisin-descent-data}, $V = T[p^{-1}]$ is a crystalline representation of $\mathcal{G}_R$. This completes the proof of Theorem~\ref{thm:main-intro}.

\bibliographystyle{amsalpha}
\bibliography{library}
	
\end{document}